\documentclass[10pt,journal,compsoc]{IEEEtran}
\usepackage{cite}

\usepackage{makeidx}  % allows for indexgeneration
\usepackage{color}
\usepackage{graphicx}

\usepackage{amsmath}
\usepackage{amssymb}
\usepackage[ruled,lined,linesnumbered]{algorithm2e}
\usepackage{algorithmic}
\usepackage{blindtext}
\usepackage{enumitem}
\usepackage{multirow}
\usepackage[labelformat=simple]{subfig}
\usepackage[pagebackref=true,breaklinks=true,letterpaper=true,colorlinks,bookmarks=false]{hyperref}
\usepackage{threeparttable/threeparttable}
\usepackage{bm}
\usepackage[normalem]{ulem}

\usepackage{amsmath}
\usepackage{ftnxtra}
\usepackage{hyperref}

\usepackage{amsfonts}
\usepackage{graphicx}
\usepackage{float}
\usepackage{graphicx}
\usepackage{epstopdf}
\usepackage{amssymb}
\usepackage{multirow}
\usepackage{amsmath}
\usepackage{algorithmic}

\newtheorem{prop}{Proposition}

\newtheorem{cor}{Corollary}
\newtheorem{exmp}{Example}
\newtheorem{rem}{Remark}
\newtheorem{pf}{Proof}
\newtheorem{df}{Definition}

\makeatletter
    
    \newcommand{\Rmnum}[1]{\expandafter\@slowromancap\romannumeral #1@}
  \makeatother
\begin{document}

%\title{\small Active Fine-Tuning: Integrating Active Learning and Transfer Learning to Reduce Annotation Efforts}
%\title{Integrating Active Learning and Transfer Learning for Biomedical Image Analysis}
\title{Frequency of Rational Fractions on [0, 1]}

\author{Zongwei~Zhou,
        Dawei Lu
\IEEEcompsocitemizethanks{
\IEEEcompsocthanksitem\normalfont Z. Zhou is with the Department of Biomedical Informatics, Arizona State University, Scottsdale, AZ 85259 USA. (\textit{zongweiz@asu.edu})
\IEEEcompsocthanksitem\normalfont D. Lu is with the School of Mathematical Sciences, Dalian University of Technology, Dalian 116023, China. (\textit{ludawei$_{-}$dlut@163.com})}
}

%\thanks{Manuscript received December, 2017.}

%\markboth{IEEE TRANSACTIONS ON PATTERN ANALYSIS AND MACHINE INTELLIGENCE,~Vol.~XX, No.~XX, XX~20XX}%
%{Shell \MakeLowercase{\textit{et al.}}: Bare Demo of IEEEtran.cls for Computer Society Journals}

\IEEEtitleabstractindextext{%
\begin{abstract}
In this paper, the authors design a trial to count rational ratios on the interval [0, 1], and plot a normalized frequency statistical graph. Patterns, symmetry and co-linear properties reflected in the graph are confirmed. The main objective is to present a new view of Farey sequence and to explain the inner principle of its procedure. In addition, we compare Farey sequence and Continued fraction in terms of numerical approximation track and clarify the internal reason why we iteratively choose mediant as the next suitable approximation for the first time. Besides, all sorts of Fibonacci-Lucas sequences emerge from the statistical graph.
\end{abstract}

\begin{IEEEkeywords}
rational fraction, normalization, frequency statistical graph, numerical approximation, Farey sequence, Continued fraction, Fibonacci-Lucas sequence
\end{IEEEkeywords}}

\maketitle

\IEEEdisplaynontitleabstractindextext

% \fi
\IEEEpeerreviewmaketitle

\section{Introduction}\label{sec1}

\IEEEPARstart{F}{or} convenience, let $\mathcal{N}$ be the set of all natural numbers, $\mathcal{N}^+$ be the set of all positive integers, $\left\lfloor x\right\rfloor$ be the floor of $x$, i.e. the largest integer less than or equal to $x$.

\begin{figure*}
\centering
    \includegraphics[width=1\textwidth]{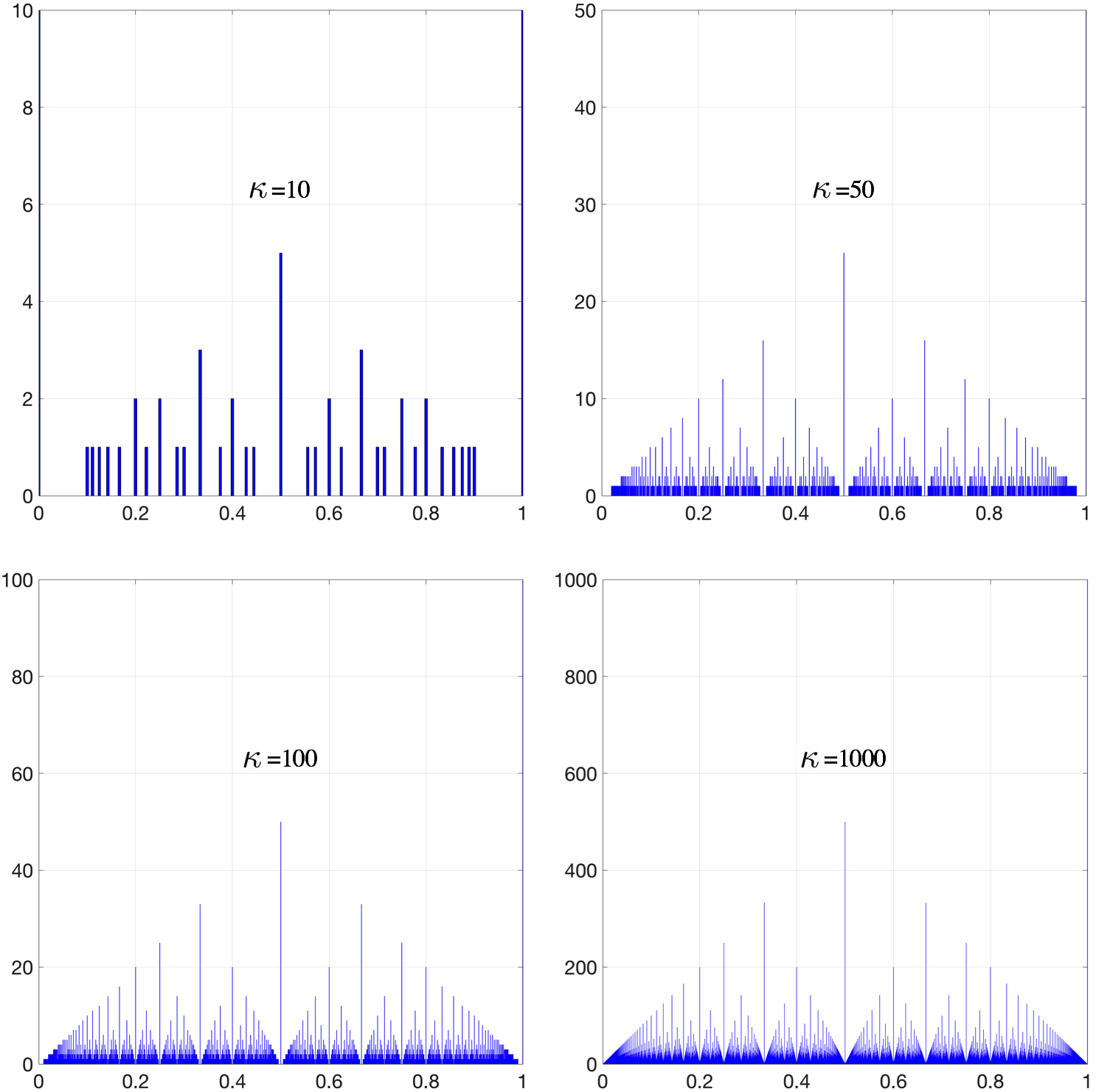}
    \caption{$\kappa$=10, 50, 100, 1000, the occurrence of each different fractions}\label{fig1}
\end{figure*}

To define a trial, where $p$ and $q$ are co-prime $(p\leq q)$, the simplified ratio form of natural numbers $a$ and $b$,
\begin{equation}\label{equ0}
\begin{split}
    T_{\kappa}(p,q)
    = \# \{ (a,b) \in \mathcal{N} \times \mathcal{N}^+ \mid a \le b \le \kappa, \\
    (a,b) \mbox{ and}  (p,q) \mbox{ are Q-scalar multiples of another.} \}
\end{split}
\end{equation}
where $\#$ represents the cardinality of a set, and $\kappa$ represents the maximum value of $q$. If we regard $(p,q)$ as a rational fraction, i.e. $p$ representing numerator and $q$ representing denominator, we will get a bunch of rational fraction sequence: $0/1$, $1/1$, $0/2$, $1/2$, $2/2$, $0/3$, $1/3$, $2/3$, $3/3$, $\cdots$, $0/\kappa$, $1/\kappa$, $2/\kappa$, $\cdots$, $\kappa/\kappa$. Then we convert these rational fractions into simplified ratios and count the occurrence times of each simplified ratio ($p$,$q$). Drawing a bar chart with X axis representing the different value of simplified fractions on the interval [0, 1], and Y axis representing individually occurrence times. As shown in Figure \ref{fig1}, we set $\kappa$ valuing 10, 50, 100, 1000.

Any number which can be expressed in the form of a quotient of two co-prime integers is called a rational number \cite{rosen2007discrete}. When $\kappa$ tends to infinity, the produced rational fractions can cover all the rational numbers on $[0, 1]$. From Figure \ref{fig1}, we observe that no matter what value $\kappa$ takes, the structure and the appearance in occurrence (Y axis) generally remains constant. When $\kappa$ increases, the rational fractions (X axis) become much more dense. To some extent, the statistical graph shows self-similarity, although not in the normal way. A self-similar object is exactly or approximately similar to a part of itself (i.e. the whole graph has the same shape as one or more of the parts) \cite{mandelbrot1967long}. From Figure \ref{fig1}, parts of them showing the asymmetric statistical properties at many scales, we named it \textbf{self-asymmetric-similarity}.

\begin{prop}\label{prop1}
    For any two co-prime integers $p$ and $q$ $(p\le q)$, we have
    \begin{equation}\label{equ1}
        T_{\kappa}(p,q)=\left\lfloor \frac{\kappa}{q}\right\rfloor.
    \end{equation}
\end{prop}

\begin{pf}
    According to the definition (\ref{equ0}), we have
    $$
        T_{\kappa}(p,q)=\# \left \{\frac{p}{q},\frac{2p}{2q},\frac{3p}{3q},\dots,\frac{\mathcal{X}p}{\mathcal{X}q}\right \}=\mathcal{X},\quad q\le \kappa.
    $$
    Since $\mathcal{X}$ is a certain integer fitting that $\mathcal{X}q\le \kappa$ and $(\mathcal{X}+1)q> \kappa$, we obtain that
    $$
        \frac{\kappa}{q}-1<\mathcal{X}\le\frac{\kappa}{q},\mbox{ thus } \mathcal{X}=\left\lfloor \frac{\kappa}{q}\right\rfloor.
    $$
The proof of Proposition \ref{prop1} is completed.
\end{pf}

\begin{algorithm}[t]
    \caption{The Procedure of Producing Fractions}
    \label{alg1}
    \begin{algorithmic}
        \REQUIRE $\kappa$
        \ENSURE each $p/q$
        \FOR{each $q\in [1,\kappa]$}
        \FOR{each $p\in [0,q]$}
        \PRINT{$p/q$}
        \ENDFOR
        \ENDFOR
    \end{algorithmic}
\end{algorithm}

Therefore, if we only account the occurrence of simplified fractions $(p,q)$, we observe that it would decrease as the denominator $q$ increases. Moreover, the symmetry is clear from Figure \ref{fig1} whatever value $\kappa$ takes. The occurrence statistical graph is symmetrical about line $x=1/2$, because rational fraction $p/q$ and $(q-p)/q$ would emerge in pairs. Namely,
\begin{equation}\label{equ2}
    T_{\kappa}(p,q)=T_{\kappa}(q-p,q).
\end{equation}
\begin{rem}

According to Proposition \ref{prop1}, $T_{\kappa}(p,q)=T_{\kappa}(q-p,q)=1/q$. Here, we only need to prove that $q-p$ and $q$ are co-prime where $p$ and $q$ are co-prime. If $(q-p,q)$ is not in the simplest form, then $\exists \lambda \ge 2$, such that $q-p=\lambda s$, $q=\lambda t$, where $s$ and $t$ are co-prime. Then we have
$$
    \frac{p}{q}=\frac{(t-s)\lambda}{t\lambda}=\frac{t-s}{t},
$$
which revealing that $p$ and $q$ are not co-prime. The assumption does not hold, and $q-p$ and $q$ are in the simplest form.
\end{rem}

The rest of the paper is arranged as follows. In Section 2, occurrence normalization is proposed and some relative properties are provided. In Section 3, all sorts of co-linear properties are classified and presented using computer assistant and statistical graph. In Section 4, the comparison between two efficient approximation methods and the inner principle of Farey procedure are given by algorithm track. Finally, in Section 5, some further works are indicated.

\section{Occurrence normalization}

We take different ratio frequencies $P_\kappa(p,q)$ into account, which equals to the occurrence time of one rational fraction $(p,q)$ divided by the total occurrence of all produced rational fractions.
\begin{equation}\label{equ3}
    P_{\kappa}(p,q)=\frac{T_{\kappa}(p,q)}{2+3+4+\dots+(\kappa+1)}=\frac{2T_{\kappa}(p,q)}{\kappa(\kappa+3)}.
\end{equation}

Different values of $P_{\kappa}(p,q)$ are shown in Table 1 with changing $\kappa$, where $(p,q)$ represents each produced simplified ratios, and $P_{\kappa}(p,q)$ represents percentage frequency of $(p,q)$. Furthermore, letting $F(\kappa)$ refers to the total number of different produced simplified ratios, $F(\kappa)$ behaves asymptotically with $F(\kappa)\sim (3\kappa^2)/\pi^2$ \cite{cobeli2003haros}.
\begin{table}[H]\caption{Different value of $P_\kappa(p,q)$ with changing $\kappa$}\label{tab1}
\begin{center}
    \begin{tabular}{|ccc|ccc|}
        \hline
            $\kappa$ & $(p,q)$ & $P_{\kappa}\times10^{-2}$ & $\kappa$ & $(p,q)$ & $P_{\kappa}\times10^{-2}$  \\
        \hline
            \multirow{4}{*}{10} & (1,1) & 15.3846 &\multirow{4}{*}{100} & (1,1) & 1.9417 \\
                & $(1,2)$ & 7.6923 & & $(1,2)$ & 0.9709 \\
                & $(1,3)$ & 4.6154 & & $(1,3)$ & 0.6408 \\
                & $(1,4)$ & 3.0769 & & $(1,4)$ & 0.4854 \\
        \hline
            \multirow{4}{*}{500} & (1,1) & 0.3976 &\multirow{4}{*}{1000} & (1,1) & 0.1994 \\
                & $(1,2)$ & 0.1988 & & $(1,2)$ & 0.0998 \\
                & $(1,3)$ & 0.1320 & & $(1,3)$ & 0.0664 \\
                & $(1,4)$ & 0.0994 & & $(1,4)$ & 0.0499 \\
        \hline
            \multirow{4}{*}{2000} & (1,1) & 0.0999 & \multirow{4}{*}{3000} & (1,1) & 0.0666 \\
                & $(1,2)$ & 0.0499 & & $(1,2)$ & 0.0333 \\
                & $(1,3)$ & 0.0333 & & $(1,3)$ & 0.0222 \\
                & $(1,4)$ & 0.0250 & & $(1,4)$ & 0.0167 \\
        \hline
    \end{tabular}
\end{center}
\end{table}

The ratio frequency would generally reduce with increasing value of $\kappa$. A small occurrence in fraction frequency results in a large percentage frequency while the reverse is true for high fraction frequencies. Obviously in the whole figure view, the general structure and trends are fixed no matter what value $\kappa$ takes. In order to see into the statistical graph more objectively, and to show each rational decimal more relatively, we take a normalization method \cite{aksoy2001feature} to each different simplified ratio frequency.
\begin{equation}\label{equ4}
    T'_{\kappa}(p,q)=\frac{T_{\kappa}(p,q)-\min(T_{\kappa})}{\max(T_{\kappa})-\min(T_{\kappa})},
\end{equation}
where $\max(T_{\kappa})$ is the maximum time that any one of simplified fractions on $[0, 1]$ can occur. Referring to Proposition \ref{prop1}, when $q=1$, the fraction either equals to $(1,1)$ or $(0,1)$, the maximum appearance time of simplified fraction is $\kappa$, i.e.,
$$
    \max(T_{\kappa})=T_{\kappa}(1,1)=T_{\kappa}(0,1)=\kappa.
$$

At the same time, $\min(T_{\kappa})$ is the minimum time that one simplified fraction can occur. Referring to Proposition \ref{prop1}, when $q=\kappa$, the minimum appearance of simplified fraction is $(1,\kappa)$, i.e.,
$$
    \min(T_{\kappa})=T_{\kappa}(1,\kappa)=1.
$$

The normalized expression adjusts to the equation below
\begin{equation}\label{equ5}
\begin{split}
    RNF(p,q)&=\lim_{\kappa\rightarrow\infty}T'_{\kappa}(p,q)=\lim_{\kappa\rightarrow\infty}\frac{T_{\kappa}(p,q)-\min(T_{\kappa})}{\max(T_{\kappa})-\min(T_{\kappa})}  \\ &=\lim_{\kappa\rightarrow\infty}\frac{T_{\kappa}(p,q)-1}{\kappa-1}=\lim_{\kappa\rightarrow\infty}\frac{T_{\kappa}(p,q)}{\kappa}.
\end{split}
\end{equation}

We named the frequency of the simplified fraction $(p,q)$ after being normalized as \textbf{Relative Normalized Frequency (RNF)}. Moreover,
$$
    RNF(p,q)=\lim_{\kappa\rightarrow\infty}\frac{T_{\kappa}(p,q)}{\kappa}\leq \lim_{\kappa\rightarrow\infty}\frac{\max(T_{\kappa})}{\kappa}=\lim_{\kappa\rightarrow\infty}\frac{\kappa}{\kappa}=1,
$$
$$
    RNF(p,q)=\lim_{\kappa\rightarrow\infty}\frac{T_{\kappa}(p,q)}{\kappa}\geq \lim_{\kappa\rightarrow\infty}\frac{\min(T_{\kappa})}{\kappa}=\lim_{\kappa\rightarrow\infty}\frac{1}{\kappa}=0.
$$

The value of $RNF$ is normalized by scaling on the interval $[0,1]$. Considering the limitation of computing, we set $\kappa=4000$ at this time. Sorting the produced simplified fractions ordered by $RNF$ decreasing size, the top 33 are shown in Table 2.
\begin{table*}[t]\caption{The top 33 simplified fractions ordered by $RNF$}\label{tab2}
\begin{center}
    \footnotesize
    \begin{tabular}{|ccc|ccc|ccc|}
        \hline
            (p,q) & $T_{\kappa}$ & $RNF$ & (p,q) & $T_{\kappa}$ & $RNF$ & (p,q) & $T_{\kappa}$ & $RNF$ \\
        \hline
            $(0,1)$   & 4000  & 1     & $(1,6)$ & 666 & $1/6$ & $(1,9)$  & 444 & $1/9$  \\
            $(1,1)$   & 4000  & 1     & $(5,6)$ & 666 & $1/6$ & $(2,9)$  & 444 & $1/9$  \\
            $(1,2)$   & 2000  & $1/2$ & $(1,7)$ & 571 & $1/7$ & $(4,9)$  & 444 & $1/9$  \\
            $(1,3)$   & 1333  & $1/3$ & $(2,7)$ & 571 & $1/7$ & $(5,9)$  & 444 & $1/9$  \\
            $(2,3)$   & 1333  & $1/3$ & $(3,7)$ & 571 & $1/7$ & $(7,9)$  & 444 & $1/9$  \\
            $(1,4)$   & 1000  & $1/4$ & $(4,7)$ & 571 & $1/7$ & $(8,9)$  & 444 & $1/9$  \\
            $(3,4)$   & 1000  & $1/4$ & $(5,7)$ & 571 & $1/7$ & $(1,10)$ & 400 & $1/10$ \\
            $(1,5)$   & 800   & $1/5$ & $(6,7)$ & 571 & $1/7$ & $(3,10)$ & 400 & $1/10$ \\
            $(2,5)$   & 800   & $1/5$ & $(1,8)$ & 500 & $1/8$ & $(7,10)$ & 400 & $1/10$ \\
            $(3,5)$   & 800   & $1/5$ & $(3,8)$ & 500 & $1/8$ & $(9,10)$ & 400 & $1/10$ \\
            $(4,5)$   & 800   & $1/5$ & $(5,8)$ & 500 & $1/8$ &        &     &        \\
                      &       &       & $(7,8)$ & 500 & $1/8$ &        &     &        \\
        \hline
    \end{tabular}
\end{center}
\end{table*}
\begin{figure}
\centering
    \includegraphics[width=\linewidth]{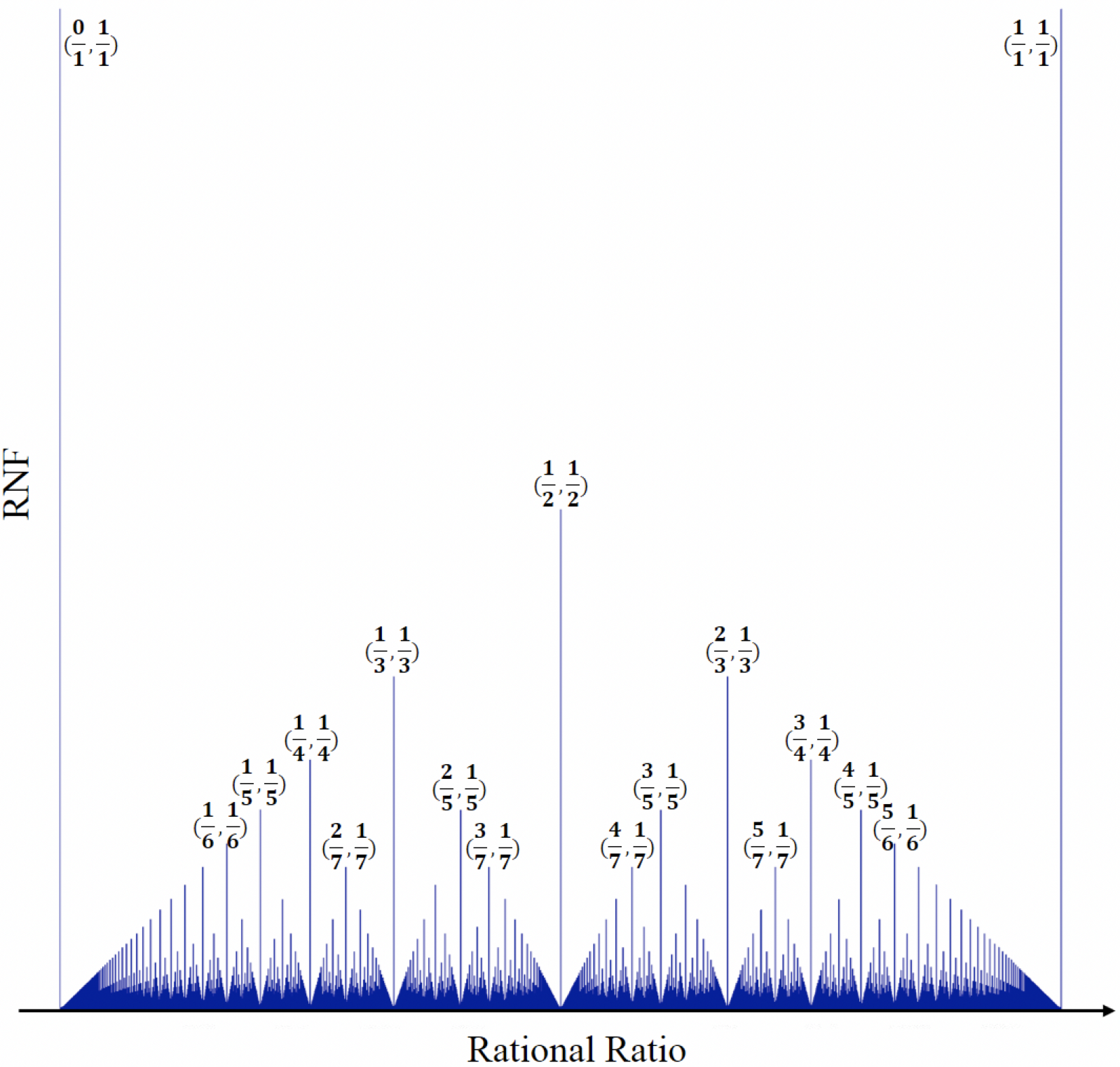}
    \caption{Normalized frequency statistical graph, $\kappa$=4000}\label{fig2}
\end{figure}

Figure \ref{fig2} is the normalized frequency statistical graph graph, X axis representing all produced simplified fractions on the interval $[0, 1]$, Y axis representing the $RNF$ corresponding to each simplified fractions appearance with the value in the interval $[0, 1]$. According to Figure \ref{fig2} and Table 2, we can draw following conclusions.
\begin{prop}\label{prop2}
    When $q$ and $p$ are co-prime integers, then
    \begin{equation}\label{equ6}
        RNF(p,q)=RNF(1,q)=\frac{1}{q},\quad q\in \mathcal{N}^+.
    \end{equation}
\end{prop}

Derived from Proposition \ref{prop1}, the larger the denominator $q$ the less chance of occurrence of the rational number $(p,q)$. Besides, Proposition \ref{prop2} extends to all the fractions between $[0,1]$, proposing that the occurrence of $(p,q)$ largely depends on denominator of fraction after simplification. Proposition \ref{prop2} can be easily proposed via Proposition \ref{prop1} and (\ref{equ5}). \\

The whole statistical graph (Figure \ref{fig2}) is symmetric about the line $x=1/2$, that is
\begin{equation}\label{equ7}
    RNF(p,q)=RNF(q-p,q),
\end{equation}
and it can be easily proposed via Proposition \ref{prop2} and (\ref{equ5}) as well.
\begin{df}\label{def1}
    Considering X coordinates, if we have two simplified fractions $x_1=(p_1,q_1)$ and $x_2=(p_2,q_2)$ in Figure \ref{fig2} with the properties that
    \begin{equation}\label{neighbor}
    |q_1p_2-p_1q_2|=1,
    \end{equation}
    all the $RNF$ of rational fractions value between $x_1$ and $x_2$ are smaller than both $x_1$ and $x_2$. Then fractions $(p_1,q_1)$ and $(p_2,q_2)$ are named as \textbf{neighbours}.
\end{df}
\begin{exmp}\label{exmp1}
    In Figure \ref{fig2}, considering X coordinates, point $1/2$ and point $1/3$ are neighbours and point $1/3$ and point $2/5$ are neighbours since none of the fractions' $RNF$ is larger than either of them. On the contrary, point $1/4$ and point $2/5$ are not neighbours because the $RNF$ of point $1/3$ is larger than these two points, and their X coordinates are not following (\ref{neighbor}).
\end{exmp}

Definition of neighbours above was firstly proposed by the British geologist John Farey, Sr., who published Farey sequences in the Philosophical Magazine in 1816 \cite{wiki2018farey}. Farey conjectured, without offering proof, that each new term in a Farey sequence expansion is the mediant of its neighbours. We use this neighbours property \cite{ainsworth2012farey} here for further exploration.
\begin{cor}\label{cor1}
    Any two neighbours have their secondary local maximum point, named as \textbf{mediant}. The mediant of these two fractions is given by
    \begin{equation}\label{mediant}
    \frac{p_1}{q_1}\oplus \frac{p_2}{q_2}=\frac{p_1+p_2}{q_1+q_2}.
    \end{equation}
\end{cor}
\begin{exmp}\label{exmp2}
     In Figure \ref{fig2}, point $2/5$ is the mediant of point $1/2$ and point $1/3$ because point $1/2$ and point $1/3$ are neighbours and point $2/5$ follows (\ref{mediant}).
\end{exmp}
\begin{pf}

Any two neighbours, with X coordinates $(p_1,q_1)$ and $(p_2,q_2)$ respectively will have their mediant $(p_1+p_2,q_1+q_2)$. Without loss of generality, we assume that $(p_1,q_1)<(p_2,q_2)$, with $p_1$ and $q_1$ co-prime, $p_2$ and $q_2$ co-prime, $p_1,p_2\in \mathcal{N}$, $q_1,q_2\in \mathcal{N}^+$.  \\

Firstly, we need to prove that
\begin{equation}\label{c1}
    \frac{p_1+p_2-1}{q_1+q_2}\leq\frac{p_1}{q_1}<\frac{p_1+p_2}{q_1+q_2}<\frac{p_2}{q_2}\leq\frac{p_1+p_2+1}{q_1+q_2},
\end{equation}
in order to demonstrate that $(p_1+p_2,q_1+q_2)$ is the only fraction whose value between $(p_1/q_1)$ and $(p_2,q_2)$ and whose denominator is no more than $q_1+q_2$.  \\

Using $(p_1,q_1)<(p_2,q_2)$, it's easy to have
$$
    \frac{p_1}{q_1}<\frac{p_1+p_2}{q_1+q_2}<\frac{p_2}{q_2}.
$$
Using $q_1p_2-p_1q_2=1$ and $q_1\geq 1$, we have
\begin{equation}\label{c2}
    \frac{p_1+p_2-1}{q_1+q_2}\leq\frac{p_1}{q_1}.
\end{equation}
Similarly, we have
\begin{equation}\label{c3}
    \frac{p_1}{q_1}\leq\frac{p_1+p_2+1}{q_1+q_2}.
\end{equation}
The proof of (\ref{c1}) is completed. \\

Secondly, we prove that $p_1+p_2$ and $q_1+q_2$ are co-prime, because only simplified fractions exists in the Figure \ref{fig2}. Presume $(p_1+p_2,q_1+q_2)$ are not in the simplest form at the first time, then $\exists \lambda \geq2$, making $p_1+p_2=\lambda s$, $q_1+q_2=\lambda t$, where $(s,t)$ are simplest form of $(p_1+p_2,q_1+q_2)$.
$$
    t=\frac{q_1+q_2}{\lambda}\leq\frac{1}{2}(q_1+q_2)\leq\max(q_1, q_2),\quad\frac{1}{t}\geq\frac{1}{\max(q_1, q_2)}.
$$
Based on Proposition \ref{prop2}, we have
$$
    RNF(p_1+p_2,q_1+q_2)=RNF(s,t)=\frac{1}{t}\geq\frac{1}{\max(q_1, q_2)}.
$$

In terms of $RNF(p_1,q_1)=(1,q_1)$ and $RNF(p_2,q_2)=(1,q_2)$, the $RNF$ of $(p_1+p_2,q_1+q_2)$ is larger than either $(p_1,q_1)$ or $(p_2,q_2)$, plus the result proved above that $(p_1+p_2,q_1+q_2)$ values between $(p_1,q_1)$ and $(p_2,q_2)$, it's contradictory to the fact that $(p_1,q_1)$ and $(p_2,q_2)$ are two neighbours. The assumption does not hold, thus $(p_1+p_2,q_1+q_2)$ are in the simplest form.  \\

Finally, we need to prove that $(p_1,q_1)$ and $(p_1+p_2,q_1+q_2)$ are neighbours, in addition, $(p_1+p_2,q_1+q_2)$ and $(p_2,q_2)$ are neighbours as well. Considering the definition of neighbours (\ref{neighbor}), we have
$$
    q_1(p_1+p_2)-(q_1+q_2)p_1=q_1p_2-q_2p_1=1,
$$
$$
    (q_1+q_2)p_2-q_2(p_1+p_2)=q_1p_2-q_2p_1=1.
$$

Talking into account all analyses above, the proof of Corollary \ref{cor1} is completed.
\end{pf}

\section{The Co-linear Points Visualization}

There are various co-linear points in Figure \ref{fig2}, especially among the local maximum points. We classify these co-linear points into two groups, named \textbf{Category A} and \textbf{Category B}. Derived from them, we can also gain the \textbf{sub-Category A} and \textbf{sub-Category B} co-linear points. Technically, all sorts of the co-linear points are tested with the help of computer assistant. This section we give specific visualization slopes of these co-linear properties.
\subsection{Category A and sub-Category A}

Taking the corresponding maximum points in Figure \ref{fig2} and plotting the respecting graphs, the result is a dispersing graph with $(0,0)$ as the point of dispersion; we classify these lines as category A. Sub-dividing the category A graph into several minimal graphs results into the sub-Category A lines. These are well illustrated in Figure \ref{fig4}.

\begin{figure*}\label{fig4}
\centering
    \includegraphics[width=\textwidth]{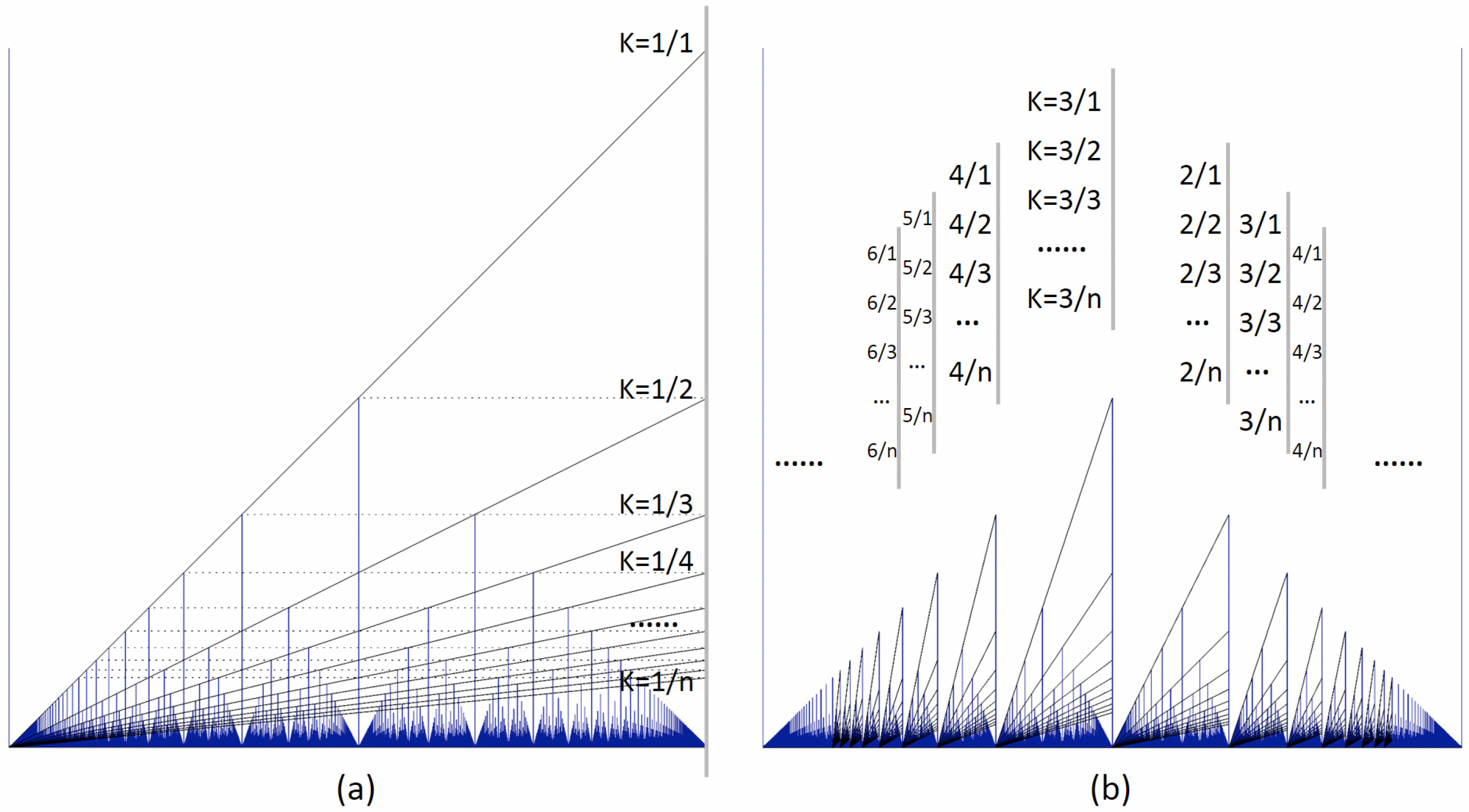}
    \caption{Co-linear points in Category A (a) and sub-Category A (b).}\label{fig4}
\end{figure*}

As seen in Figure \ref{fig4}$(a)$, the slope $k$ of the graphs is represented by $1/n$, with $n$ unit increasing from one to infinity. In Figure \ref{fig4}$(b)$, the slope of respective graphs can be defined by as the quotient of a constant numerator and a denominator of value increasing from one to infinity. The numerator value of the slopes unit increases from three on the left of the symmetrical axis and unit increases from two on the right.
\subsection{Category B and sub-Category B}

Taking the corresponding maximum points in Figure \ref{fig2} and plotting the respecting graphs, the result is a dispersing graph with $(1,1)$ as the point of dispersion; we classify these lines as category B. Sub-dividing the category B graph into several minimal graphs results into the sub-Category B lines. These are well illustrated in Figure \ref{fig5}.   \\

\begin{figure*}\label{fig5}
\centering
    \includegraphics[width=\textwidth]{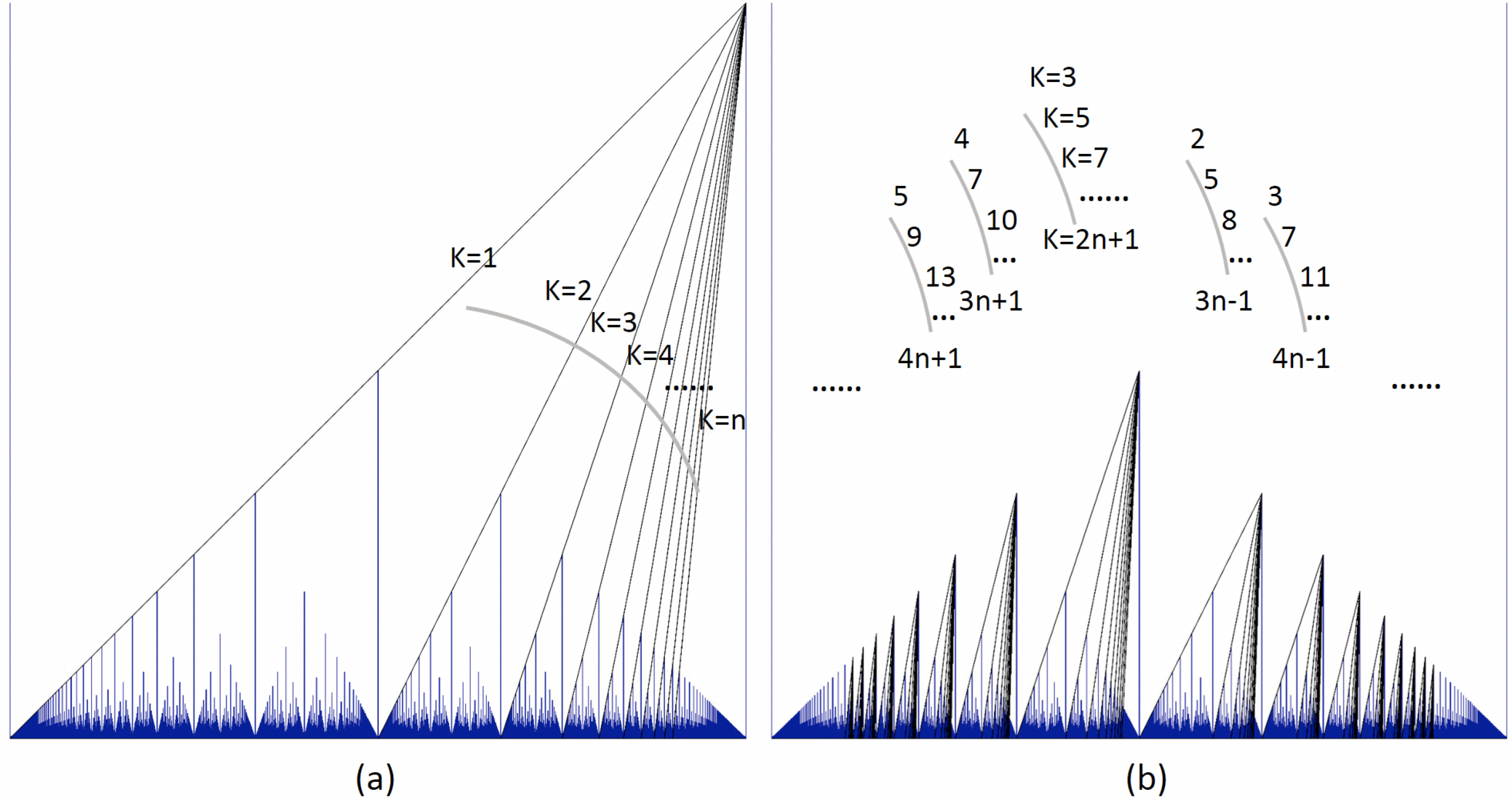}
    \caption{co-linear points in Category B (a) and sub-Category B (b)}\label{fig5}
\end{figure*}

The sub-category B graphs can be divided into left and right graphs with the line $x=1/2$ as the mirror line. We can deduce that the left graphs have a general slope of $k=an+1$ and the right graphs a slope of $k=an-1$. The co-efficient $a$ represents the respective sub-category number. For left graphs, $a(a\geq2)$ starting from the central sub-category, and $a(a\geq3)$ for right graphs with the first sub-category adjacent to the central subset. It has a unit increase as we proceed to adjacent sub-category graphs. In addition, $n$ represents the $n^{th}$ line in a sub-category set.  \\

Therefore, the second line in the central sub-category, for example, will have a slope of $k=2\times2+1=5$, while the first line in the third right sub-category will have a slope of $k=4\times1-1=3$.

\section{Numerical approximations}

Given an irrational on the interval $[0,1]$, we can find a rational approximation using two efficient strategies, i.e. Farey sequence and Continued fraction.

A Farey sequence $F_n$ is the set of rational numbers $(p,q)$ with $p$ and $q$ co-prime, with $0\leq p\leq q\leq n$, ordered by size. For instance,

\begin{equation}
    \begin{aligned}
        &F_1=\{0/1, ~~~~~~~~~~~~~~~~~~~~~~~~~~~~~~~~~~~ 1/1\},  \\
        &F_2=\{0/1, ~~~~~~~~~~~~~~~ 1/2, ~~~~~~~~~~~~~~ 1/1\},  \\
        &F_3=\{0/1, ~~~~~~~~ 1/3,~1/2,~ 2/3,~~~~~~~~1/1\},  \\
        &F_4=\{0/1,~1/4,~1/3,~1/2,~2/3,~3/4,~1/1\}.  \\
    \end{aligned}
\end{equation}

According to Farey's definition, the mediant of two neighbours is the Freshman's sum of those neighbours, which coincides to Corollary \ref{cor1}. Nevertheless, it's a simple rule without any reason why we produce mediant in this particular method. Moreover, with the Stern-Brocot Tree, which was first described by Moritz Stern in 1858 \cite{ainsworth2012farey,graham1989concrete,bates2010locating}, given an irrational number between 0 and 1, it's possible to find a rational approximation using iterations of the mediant in the Farey sequence. Specific approximation example given by Cook \cite{john2010app}.

Continued fraction are one of the most revealing representations of numbers \cite{wall1948analytic}. To define a general continued fraction of a number as
\begin{equation}\label{equ8}
    [a_1,a_2,a_3,\cdots,a_N]=1/(a_1+1/(a_2+1/(a_3+\dots)))\quad (a_i\in \mathcal{N}),
\end{equation}
when we truncate a continued fraction after some number of terms, we get what is called a convergent \cite{lamb_2015}. We denote them by $(p,q)$. As $q$ increases, the different between target and $(p,q)$ decreases.  \\
\begin{algorithm}
    \caption{Farey Sequence Method}
    \label{alg2}
    \begin{algorithmic}
        \REQUIRE $Decimal$, $error$
        \ENSURE $Loop$, $Expansion$
        \STATE $lp \gets 0$, $lq \gets 1$
        \STATE $hp \gets 1$, $hq \gets 1$
        \STATE $mp \gets lp + hp$, $mq \gets lq + hq$
        \STATE $Loop \gets 1$
        \WHILE{$|mp / mq - Decimal| \geq error$}
        \STATE $Loop \gets Loop + 1$
        \IF{$mp / mq > Decimal$}
        \STATE $hp \gets mp$, $hq \gets mq$
        \ELSE
        \STATE $lp \gets mp$, $lq \gets mq$
        \ENDIF
        \STATE $mp \gets lp + hp$, $mq \gets lq + hq$
        \STATE $Expansion \gets mp / mq$
        \PRINT $Expansion$
        \ENDWHILE
    \end{algorithmic}
\end{algorithm}
\begin{algorithm}
    \caption{Continued Fraction Method}
    \label{alg3}
    \begin{algorithmic}
        \REQUIRE $Decimal$, $error$
        \ENSURE $Loop$, $Expansion$
        \STATE $Loop \gets 1$
        \STATE $mod \gets 1$
        \STATE $rad\gets Decimal$
        \STATE $i\gets 0$
        \WHILE{$mod \geq error$}
        \STATE $contFraction[i]\gets 0$
        \WHILE{$rad-conFraction[i]*mod\geq mod$}
        \STATE $Loop\gets Loop+1$
        \STATE $conFraction[i]\gets conFraction[i]+1$
        \STATE $Expansion\Leftarrow conFraction$
        \PRINT $Expansion$
        \ENDWHILE
        \STATE $temp\gets mod$
        \STATE $mod \gets rad-conFraction[i]*mod$
        \STATE $rad\gets temp$
        \STATE $i\gets i+1$
        \ENDWHILE
    \end{algorithmic}
\end{algorithm}

Track the expansion fractions of these two method, Farey sequence and Continued fraction, using computer programming with Algorithm \ref{alg2} and Algorithm \ref{alg3}. Despite the different definitions, their numerical approximation processes and the time complexities of the algorithms are exactly the same. For example, in the case of the numerical approximations to $(\sqrt{5}-1)/2$ provided by Farey sequence, the track is $(1,2)$, $(2,3)$, $(3,5)$, $(5,8)$, $(8,13)$, $(13,21)$, ... iterating through the Stern-Brocot Tree. Meanwhile, the track of continued fraction method is $[1]$, $[1,1]$, $[1,1,1]$, $[1,1,1,1]$, $[1,1,1,1,1]$, $[1,1,1,1,1,1]$, ... Transformed into simplified rational fraction form, both of them meet one certain approximation track. According to the testing result, these two approximation procedures are intrinsic isogenesis.

The question is why we choose the mediant fraction as the more suitable approximation. These algorithms have wide applications \cite{rat2017matlab}, though, we haven't realized the inner mechanisation of them and why they work so well. Taking into consideration the procedure of producing fractions in (\ref{equ0}), we merely list all the natural numbers in a simple rule then count each value of fractions' occurrences. The statistical graph turns out to present such a similar outline with Farey sequence and moreover, the mediant between two neighbours is exactly following the rule of the mediant procedure as Farey proposed. It is the fact that the frequency of mediant is the local maximum point between its neighbours that leads us to believe that if we choose the mediant, it has more chances to hit the target. And even if it isn't the target, we would choose the maximum frequency point on a much smaller interval.

\section{Further works}

The difference between the mathematical thinking and the computing thinking is that we have got to regard the material limitation and the space-time limitation of the calculations. In computer language, calculations always have a definite value whereas in a mathematical language, we could limit these calculations to an infinite concept. This is the beauty of the combination of mathematical derivation and computer programming. In this section, we list some potential orientations of further works.

\subsection{Fibonacci sequence in statistical graph}

The Fibonacci sequence, also known as the Golden series, is defined recursively as following: $F(0)=0$, $F(1)=1$, $F(n)=F(n-1)+F(n-2)$, $n\geq2$, $n\in \mathcal{N}^+$, discovered by an Italy mathematician Leonardo Fibonacci \cite{azarian1990generating}.  \\
\begin{figure}
\centering
    \includegraphics[width=\linewidth]{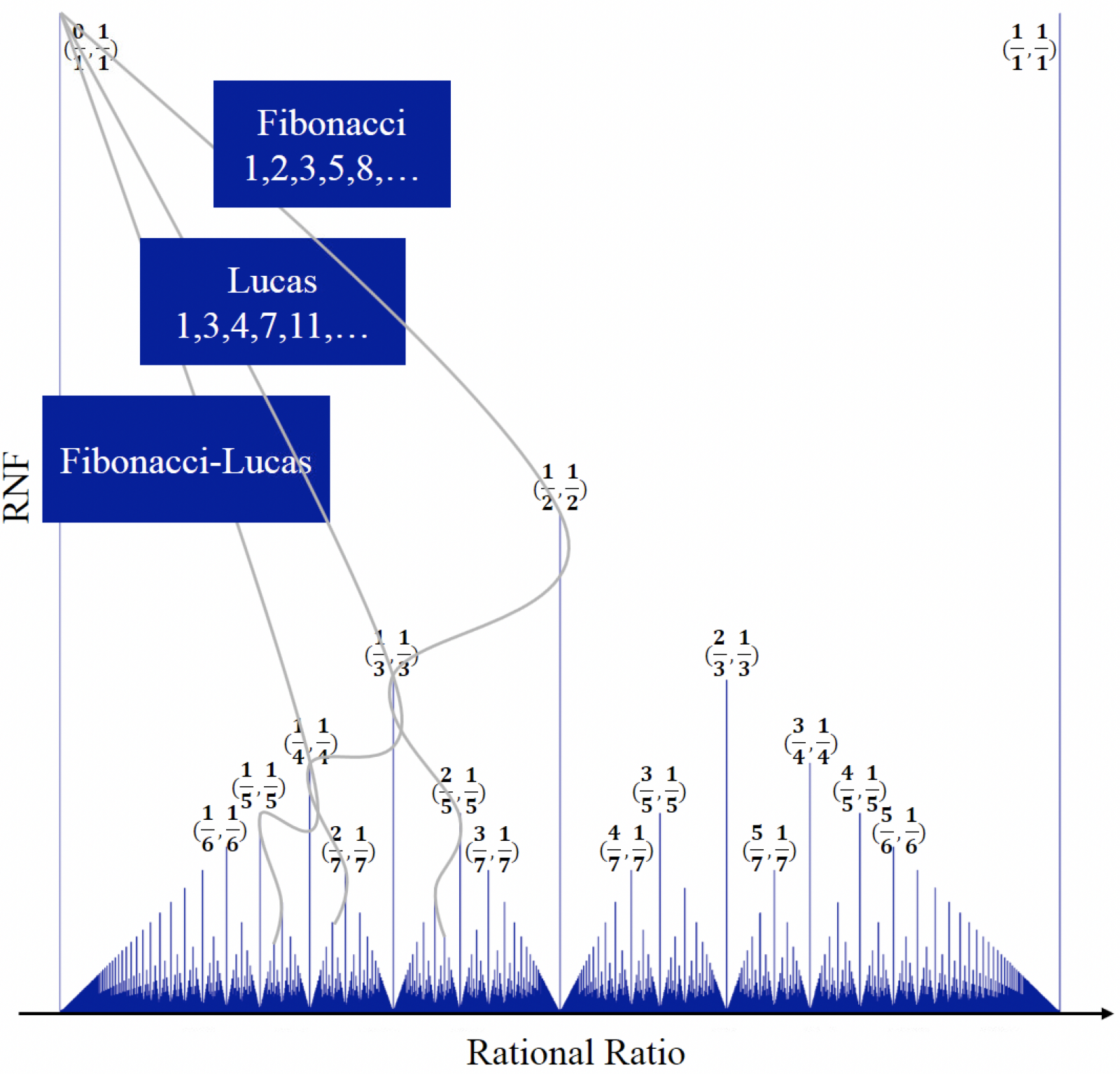}
    \caption{Fibonacci sequence in statistical graph}\label{fig6}
\end{figure}

In Figure \ref{fig2}, each mediant ratio has recursively relationship with its neighbors, i.e. mediant numerator equals to neighbors numerators added up and mediant denominator equals to neighbors denominators added up according to Corollary \ref{cor1}. There are sorts of similar series, uniformly called the Fibonacci-Lucas sequence \cite{vajda1989fibonacci}. Focusing all the numerators and the denominators in Figure \ref{fig6}, they belong to Fibonacci-Lucas sequence, and they are demonstrated to us all perfectly. Each elements in the Fibonacci-Lucas sequence show in the $\mathcal{S}$ pattern trend. What's more, different type of $\mathcal{S}$ can name its corresponding Fibonacci-Lucas sequence. According to the symmetry, there is the same rule at the right part of line $x=1/2$. It's another way to present the beauty of the statistical graph. We could infer that even some Golden phenomenon exists in the statistical graph.

\subsection{What's the other side like}

\begin{figure}
\centering
    \includegraphics[width=\linewidth]{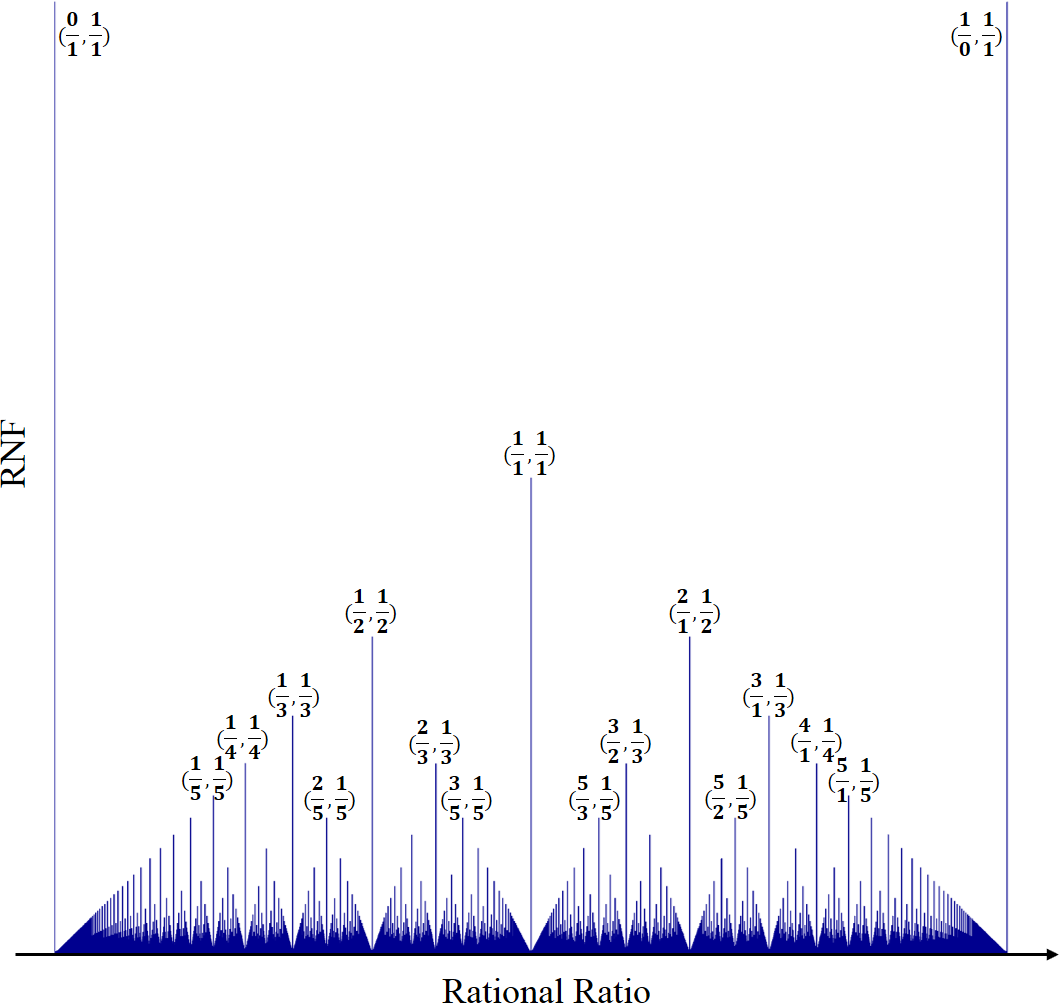}
    \caption{Quotient Statistical graph expand to positive numbers}\label{fig7}
\end{figure}

During the process of research above, we only take the rational fraction between $[0, 1]$ into account. If it is extended to an arbitrary non-negative rational decimal, i.e.,

\begin{equation}
\begin{split}
    T_{\kappa}(p,q)
    = \# \{ (a,b) \in \mathcal{N} \times \mathcal{N}^+ \mid a \le \kappa, b \le \kappa, \\
    (a,b) \mbox{ and}  (p,q) \mbox{ are Q-scalar multiples of another.}, \}
\end{split}
\end{equation}
the sequence extends to become $0/1$, $1/1$, $2/1$, $3/1$, $\cdots$, $\kappa/1$, $0/2$, $1/2$, $\cdots$, $\kappa/2$, $\cdots$, $\kappa / \kappa$, in which $\kappa$ is a huge integer tends to infinity. When $\kappa$ approaches infinity, the quotient can theoretically cover all of the non-negative rational decimal. If we recorded all of these rational ratios frequency of occurrence, what it would be? Luckily, we use computer aided, and set $\kappa$ equals to 1000, after normalized method, sorting the rational fractions according to their frequency of occurrence, then we take the top 36, and get Figure \ref{fig7}.   \\

X axis representing all produced positive simplified fractions, Y axis representing the $RNF$ corresponding to each rational decimal appearance with the value between 0 and 1. We find that the sub-graph among $[1, 2]$, $[2, 3]$, $[3, 4]$ and $[n, n+1]$ have a self-asymmetric-similarity to statistical graph in $[0, 1]$; in the latter the statistical graph is highly symmetric while in the former the statistical graph is approximately symmetric. We can draw the following conclusion with the help of computer assisting that all the Propositions and Corollaries are fit in this larger statistical graph in [0, 1]. Furthermore, we could expect that the final statistical graph that embrace all the simplified fractions including positive, negative and zero are, as double of the Figure \ref{fig7}, symmetric with line x=0. At the end of this part, we list some propositions that are certified by computer programming.
\begin{prop}\label{prop5}
    An rational number between [1, $\infty$] has its own counterpart rational number between [0, 1], whose X coordinates are reciprocal, and y coordinates are the same. For any two co-prime integers $p$ and $q$ $(p\le q)$, we have
    \begin{equation}
        \left(\frac{p}{q},\frac{1}{q}\right) \leftrightarrow \left(\frac{q}{p},\frac{1}{q}\right).
    \end{equation}
\end{prop}
\begin{prop}\label{prop6}
    In [0, 1], each maximum point has linear properties within the statistical graph; in [1, $\infty$], each maximum point has inverse properties, and these properties are not only among the maximums but among the second maximums and so on.
\end{prop}

\section*{Acknowledgments}
The research of the first author is Supported by the National Natural Science Foundation of China (Grant Nos. 11371077 and 11571058) and the Fundamental Research Funds for the Central Universities (Grant No. DUT15LK19). Computations made in this paper were performed using Matlab software.

\bibliographystyle{IEEEtran}
\bibliography{refs}

\iffalse

\fi
\end{document}